\documentclass[11pt]{article}
\usepackage[a4paper,margin=24mm]{geometry}
\usepackage[T1]{fontenc}
\usepackage{lmodern}
\usepackage{amsmath,amssymb,amsthm,mathtools,bm,booktabs,microtype}
\usepackage[colorlinks=true,linkcolor=blue,citecolor=blue,urlcolor=blue]{hyperref}
\hypersetup{pdfauthor={Robert Vrabel},pdftitle={Spectral Exactness and Quantitative Instability Approximation in the Euler--Zeitlin Model for Zonal Equilibria}}
\usepackage{enumitem}
\allowdisplaybreaks
\newcommand{\shorttitle}{Spectral Exactness for Zonal Euler Flows}
\numberwithin{equation}{section}
\newtheorem{theorem}{Theorem}[section]
\newtheorem{proposition}[theorem]{Proposition}
\newtheorem{lemma}[theorem]{Lemma}
\newtheorem{corollary}[theorem]{Corollary}
\theoremstyle{remark}

\newcommand{\HH}{\mathcal H}
\newcommand{\DD}{\mathcal D}

\newcommand{\norm}[1]{\left\lVert #1\right\rVert}
\newcommand{\ip}[2]{\left\langle #1,#2\right\rangle}
\newcommand{\spec}{\sigma}

\newcommand{\id}{\mathrm I}
\newcommand{\Rea}{\operatorname{Re}}
\newcommand{\C}{\mathbb C}
\newcommand{\R}{\mathbb R}
\newcommand{\Z}{\mathbb Z}
\title{Spectral Exactness and Quantitative Instability Approximation\\in the Euler--Zeitlin Model for Zonal Equilibria}
\author{Robert Vrabel\\[3pt]
\small Slovak University of Technology in Bratislava\\
\small Institute of Applied Informatics, Automation and Mechatronics\\
\small Bottova 25, 917 24 Trnava, Slovakia\\
\small \texttt{robert.vrabel@stuba.sk}%\\
}
\date{}
\begin{document}
\markboth{\shorttitle}{\shorttitle}
\maketitle
\begin{abstract}
Building on finite-time convergence results for geometric matrix hydrodynamics,
we establish global spectral exactness for Zeitlin's approximation of the
two-dimensional Euler equation on the sphere linearized about polynomial zonal
vorticities. By separating each azimuthal block into a skew-adjoint transport
operator and a compact perturbation, we obtain a dimension-uniform estimate
$\|K_N-K\|=O(N^{-1})$ and the sharper fixed-block estimate
$\|K_{N,m}-K_m\|=O(N^{-2})$. These bounds, together with strong-resolvent
convergence of transport, yield exact approximation of the isolated spectrum
off the imaginary axis: algebraic multiplicities are preserved, Riesz
projections converge in norm, and spurious off-axis spectral limits are
excluded. Every algebraically simple off-axis eigenvalue has an
$O(N^{-2})$ approximation error. For eigenvalues of higher multiplicity
we establish second-order convergence of arithmetic means of eigenvalue clusters and
$O(N^{-2/\alpha})$ individual eigenvalue bounds, where $\alpha$ is
the maximal Jordan-block size. For any prescribed positive growth margin,
the numbers of unstable and stable eigenvalues beyond that margin are
eventually correct. The argument combines an exact matrix-harmonic coordinate
recurrence, uniform bounds for high azimuthal modes, and resolvent factorizations
involving compact operators. An explicit degree-three zonal equilibrium has a rigorously
enclosed simple unstable eigenvalue, providing a concrete instance of the
quantitative theory.
\end{abstract}

\section{Introduction}
The two-dimensional Euler equation on the sphere combines vorticity
transport with the geometry of area-preserving flows. Zeitlin's matrix
hydrodynamics replaces the infinite-dimensional Poisson algebra by
finite-dimensional matrix commutators while retaining the Lie--Poisson
structure at each resolution~\cite{Zeitlin1991,ModinPerrot,ModinViviani2026}.
This makes the model particularly attractive for studying the onset of
instability in spherical flows. Yet geometric structure preservation and
convergence of trajectories do not, by themselves, establish spectral
fidelity: an eigenvalue computed from a finite matrix linearization may
signal exponential growth even when no corresponding mode exists in the
continuous Euler linearization. The central question is therefore whether
an observed matrix instability is a genuine Euler instability, with the
correct growth rate, multiplicity, and invariant spectral subspace.

The distinction is substantive for the Euler equation. Modin and
Perrot~\cite{ModinPerrot} proved convergence of Eulerian and Lagrangian
Jacobi perturbations on bounded time intervals, while Modin and
Viviani~\cite{ModinViviani2026} developed further convergence results for
matrix hydrodynamics. The linearized transport generator, however, is
unbounded on the full azimuthal direct sum, and its essential spectrum is
confined to the imaginary axis. Moreover, the linearization need not be
normal, so small residuals and accurate finite-time trajectories alone
do not control isolated eigenvalues or algebraic multiplicities. A
spectral approximation theorem must exclude spurious eigenvalue limits
and retain spectral subspaces in an operator topology strong enough to
handle nonnormality.

We establish this spectral correspondence for real polynomial zonal
vorticities on $\mathbb S^2$. Every such profile is a stationary
Euler flow and quantizes to an exact stationary Euler--Zeitlin flow.
Their linearizations admit a common harmonic representation and an
azimuthal direct-sum decomposition
\[
L=\bigoplus_{m\in\mathbb Z}(A_m+K_m),\qquad
L_N=\bigoplus_{m\in\mathbb Z}(A_{N,m}+K_{N,m}),
\]
where $A_m$ and $A_{N,m}$ are skew-adjoint and $K_m$ and $K_{N,m}$ are
compact or finite-rank perturbations involving the inverse Laplacian.
The main analytical obstacle is uniformity across this direct sum.
Fixed-azimuthal consistency does not control modes whose azimuthal index
grows with matrix dimension. The inverse-Laplacian structure and an exact
matrix-coordinate recurrence together provide the necessary global control.

\noindent\textbf{Principal results.}
We prove five quantitative statements for this complete direct-sum
problem. First, the compact parts converge in operator norm globally,
while every fixed nonzero azimuthal block admits a stronger second-order
rate:
\begin{equation}\label{eq:intro-rates}
\|K_N-K\|_{\mathcal L(\HH)}\le \frac{C}{N},\qquad
\|K_{N,m}-K_m\|_{\mathcal L(\HH_m)}\le \frac{C_m}{N^2}.
\end{equation}
Second, the Zeitlin spectra are \emph{exact away from the imaginary
axis}: each isolated Euler eigenvalue is approximated with its full
algebraic multiplicity, the corresponding Riesz projections converge in
operator norm, and every convergent sequence of off-axis matrix
eigenvalues has an Euler spectral limit. Third, a simple isolated
eigenvalue $\lambda$ has a matching matrix eigenvalue satisfying
$|\lambda_N-\lambda|\le C_\lambda N^{-2}$. Fourth, an isolated eigenvalue
of algebraic multiplicity $\mu$ and maximal Jordan-block size $\alpha$
has a matrix cluster whose arithmetic mean converges at order $N^{-2}$,
while its individual eigenvalues converge at order $N^{-2/\alpha}$;
the generalized eigenspace is approximated at second order on the
limiting spectral subspace. Fifth, for any fixed growth threshold
$\delta>0$ not crossed by the Euler spectrum, the finite matrix model
eventually has exactly the correct number of eigenvalues in each of
the half-planes $\Re z>\delta$ and $\Re z<-\delta$, counted
algebraically. Thus the matrix approximation recovers both the \emph{location}
and the \emph{total algebraic multiplicity} of exponentially unstable modes beyond any
fixed positive growth margin. These are properties of the entire matrix
spectrum, not merely of a preselected Fourier block.

The proof has two complementary components. On the algebraic side,
exact spin-matrix recurrences and trace-orthogonal zonal polynomials
give the coordinate-consistency bounds, including all top-degree
truncation terms. An inverse-Laplacian factor yields uniform decay of
the compact perturbation at large azimuthal number and thereby
promotes fixed-block estimates to global operator-norm convergence.
On the spectral side, strong-resolvent convergence of the skew-adjoint
transport parts combines with norm convergence of compact
Birman--Schwinger operators. A compact correction to the full resolvent
then establishes norm convergence of Riesz projections. The
quantitative eigenvalue estimates follow from smoothness of isolated
off-axis root spaces, consistent blockwise discretization, and
finite-dimensional perturbation theory for Jordan matrices.

The setting is connected to established analyses of zonal stability on
the sphere, including Taylor's degree-three Jacobi-operator analysis,
with an appendix by Marzuola and Taylor~\cite{Taylor2016}, and to spectral studies
of structure-preserving truncations for Euler flows on the
torus~\cite{Dullin2016}. Contemporary analyses of rotating-sphere zonal
flows emphasize the relevance of low-degree jets and their instability
thresholds~\cite{ConstantinGermain2022,ConstantinEtAl2025}.
The present results address the unforced Euler--Zeitlin approximation,
and establish a quantitative operator-theoretic connection between its
finite-dimensional spectra and the spherical Euler linearization.
The spectral reduction is related to classical compact-operator approximation
theory~\cite{Kato,Osborn1975,Anselone}, while the uniform quantized-zonal estimates
are specific to the matrix hydrodynamics considered here. The present estimates concern the spherical
matrix quantization, the entire azimuthal direct sum, and its isolated
off-axis spectrum. We conclude with a certified unstable eigenvalue for
a degree-three profile, which links the quantitative approximation theory
to a concrete spectral instability.

The paper is organized as follows. Section~\ref{sec:setting} fixes the
Euler--matrix identification. Section~\ref{sec:poly} develops the exact
zonal and coordinate algebra. Section~\ref{sec:bounds} proves the global
and blockwise norm estimates. Sections~\ref{sec:spectral} and
\ref{sec:eig} establish spectral exactness and quantitative
approximation of simple eigenvalues and multiple spectral clusters. Section~\ref{sec:example} gives the unstable-mode
enclosure, Section~\ref{sec:discussion} draws the main consequences,
and the rational bounds are specified in the appendix.

\section{Euler and Zeitlin linearizations on a common space}\label{sec:setting}
Let $\mathbb S^2$ be the unit sphere with standard area element and let $Y_{\ell m}$ be complex spherical harmonics orthonormal in $L^2(\mathbb S^2)$. We use the Laplacian convention
\begin{equation}\label{eq:lap}
-\Delta Y_{\ell m}=\ell(\ell+1)Y_{\ell m},\qquad \ell\ge1,\quad |m|\le\ell.
\end{equation}
Set $\HH=L^2_0(\mathbb S^2;\C)$ and
\begin{equation}\label{eq:blocks}
\HH=\bigoplus_{m\in\Z}\HH_m,\qquad
\HH_m=\overline{\operatorname{span}}\{Y_{\ell,-m}:\ell\ge\max(1,|m|)\}.
\end{equation}
We use the spectral Sobolev norm $\|u\|_{H^s}^2=\sum_{\ell\ge1}\sum_{|m|\le\ell}(1+\ell(\ell+1))^s|u_{\ell m}|^2$. Throughout, let
\begin{equation}\label{eq:base}
\omega_*(x)=\sum_{\ell=1}^{d}b_\ell Y_{\ell0}(x),\qquad
x=\cos\theta,\qquad b_\ell\in\R,\quad d\text{ fixed},\qquad
\psi_*=\Delta^{-1}\omega_*.
\end{equation}
All derivatives of $\omega_*$ and $\psi_*$ below refer to the scalar coordinate $x$. Since both depend only on $x$, $\{\psi_*,\omega_*\}=0$. Thus $\omega_*$ is stationary.

We adopt the physical Poisson bracket of Modin--Perrot,
\begin{equation}\label{eq:physicalbracket}
\{f,g\}=\mathbf x\cdot(\nabla f\times\nabla g),\qquad
\{f(x),u(x)e^{ik\varphi}\}=-ikf'(x)u(x)e^{ik\varphi}.
\end{equation}
The physical Fourier number is $k$. In the decomposition~\eqref{eq:blocks}
we have deliberately used the reversed index $m=-k$; write
$Z_{\ell m}:=Y_{\ell,-m}$. This unitary relabeling matches the
matrix-unit index $m=r-s$. Consequently,
linearization of $\partial_t\omega+\{\psi,\omega\}=0$,
$\Delta\psi=\omega$, has the representation
\begin{equation}\label{eq:lm}
L_m=A_m+K_m,\qquad
A_m=-im\psi_*'(B_m),\qquad
K_m=-im\omega_*'(B_m)D_m^{-1},\quad m\ne0,
\end{equation}
where $B_m$ is multiplication by $x$ and
\begin{equation}\label{eq:D}
D_mZ_{\ell m}=\ell(\ell+1)Z_{\ell m},\qquad
\norm{D_m^{-1}}=\frac1{|m|(|m|+1)}\quad(m\ne0).
\end{equation}
The Fourier reindexing in~\eqref{eq:blocks} fixes the orientation consistently in both models. On $\HH_0$ the linearization vanishes. Each $A_m$ is bounded and skew-adjoint on $\HH_m$; the direct sum $A=\bigoplus_m A_m$ is generally unbounded with its natural domain
\begin{equation}\label{eq:Adomain}
\DD(A)=\bigl\{u=(u_m)\in\HH:\sum_m\norm{A_mu_m}^2<\infty\bigr\}.
\end{equation}
The operator $A$ is skew-adjoint, and $K=\bigoplus_mK_m$ is bounded by Lemma~\ref{lem:highm} below.

For $N\ge2$ put $j=(N-1)/2$, $J=j(j+1)=(N^2-1)/4$, and $\hbar_N=J^{-1/2}$. Let $p_N$ denote the skew-Hermitian matrix-harmonic quantization, with
$p_NY_{\ell m}=iT^{(N)}_{\ell m}$. For spectral comparison we use the
explicit unitary identification $\mathcal U_NT^{(N)}_{\ell m}=Z_{\ell m}
=Y_{\ell,-m}$, $1\le\ell<N$, so the matrix weight $m=r-s$ corresponds to
physical Fourier number $-m$. This is the unitary harmonic identification with reversed azimuthal index; it is distinguished from the canonical quantization map $p_N$.
The factor $i$ is separated from the Hermitian zonal harmonics below. The matrix Laplacian has eigenvalue $-\ell(\ell+1)$ on $T^{(N)}_{\ell m}$ \cite{ModinPerrot}. We henceforth write
\begin{equation}\label{eq:PN}
P_N:\HH\to\HH,\qquad P_Nu=\sum_{\ell=1}^{N-1}\sum_{m=-\ell}^{\ell}\ip{u}{Y_{\ell m}}Y_{\ell m}.
\end{equation}
The quantized equilibrium $W_{*,N}=p_N\omega_*$, $N>d$, is diagonal in the $J_3$ basis; the quantization of $\psi_*$ is diagonal as well and the two matrices commute. Consequently $W_{*,N}$ is exactly stationary.

We denote by $L_N=A_N+K_N$ the matrix linearization, transported by $\mathcal U_N$ into $P_N\HH$ and extended by zero on $(I-P_N)\HH$. Its $m$-blocks are $A_{N,m}$ and $K_{N,m}$. Both vanish for $|m|\ge N$, and $A_N^*=-A_N$. The action of $K_{N,m}$ is a commutator with $W_{*,N}$ followed by the inverse matrix Laplacian; the order of these factors is essential.

We systematically distinguish three maps: the infinite-dimensional coordinate multiplier $B_m$, the finite anticommutator coordinate $B_{N,m}$, and the harmonic projection $P_{N,m}$. In particular, $B_{N,m}$ is defined \emph{only} for $m\ne0$, where it preserves the corresponding off-diagonal matrix-unit space. All matrix norms below use the normalized trace/Hilbert--Schmidt inner product; the harmonic identification is unitary for this inner product. The discrete Laplacian agrees with $-D_m$ on $P_{N,m}\HH_m$.

\begin{lemma}[Euler--matrix orientation and block identification]\label{lem:orientation}
With the bracket~\eqref{eq:physicalbracket}, the physical Euler block at
Fourier number $k$ is
\begin{equation}\label{eq:physL}
\widetilde L_k=ik\bigl(\psi_*'(B_k)+\omega_*'(B_k)D_k^{-1}\bigr).
\end{equation}
The $m$th matrix-weight block ($r-s=m$), when transported by $\mathcal U_N$,
has the signed centered-difference representation~\eqref{eq:exactfinite}.
The latter is consistent with the physical block $\widetilde L_{-m}$.
For the solid-body rotation check $\omega_*(x)=cx$, $c=\sqrt{3/(4\pi)}$,
the physical $k=1$ block satisfies
$\widetilde L_1Y_{21}=-icY_{21}/3$, and the $m=-1$
matrix-weight block has the same limiting value on its degree-two mode.
\end{lemma}
\begin{proof}
The identity $\{f(x),u(x)e^{ik\varphi}\}=-ikf'(x)u e^{ik\varphi}$
follows from $\partial_\theta f=-\sin\theta f'(x)$ and
$\mathbf x\cdot(\nabla f\times\nabla g)
=\sin^{-1}\theta(\partial_\theta f\,\partial_\varphi g
-\partial_\varphi f\,\partial_\theta g)$.
For a vorticity perturbation $v$, the stream perturbation is
$\delta\psi=-D_k^{-1}v$. Thus
$-\{\psi_*,v\}-\{\delta\psi,\omega_*\}
=ik(\psi_*'v+\omega_*'D_k^{-1}v)$,
which gives~\eqref{eq:physL} and then~\eqref{eq:lm} upon setting $k=-m$.

For the matrix equation $\dot W+\sqrt J[P,W]=0$, $\Delta_NP=W$,
write the diagonal equilibrium as $P_*=iS_N(J_3/\sqrt J)$ and
$W_*=iQ_N(J_3/\sqrt J)$. A variation $\delta W=iZ$ has
$\delta P=-iD^{-1}Z$, and its linearized equation is
\[
\dot Z=-i\sqrt J\bigl([S_N(J_3/\sqrt J),Z]
+[Q_N(J_3/\sqrt J),D^{-1}Z]\bigr).
\]
On $E_{rs}$ with $r-s=m$, both commutators multiply by the diagonal
entry difference; this yields~\eqref{eq:exactfinite} exactly, including the
factor $-i$. The coordinate recurrence in Lemma~\ref{lem:rec} uses only
$m^2$ and is unchanged when the physical harmonics are reindexed as
$Y_{\ell,-m}$. Finally, $\psi_*=-cx/2$, so
$\widetilde L_1Y_{21}=i(-c/2+c/6)Y_{21}=-icY_{21}/3$;
the matrix-weight index $m=-1$ gives the identical leading factor.
\end{proof}

\section{Zonal polynomial quantization and the coordinate recurrence}\label{sec:poly}
For fixed $\ell$, the Hermitian diagonal matrix $T^{(N)}_{\ell0}=-i p_NY_{\ell0}$ is a real polynomial in $J_3/\sqrt J$. Write its scalar polynomial as $q_{\ell,N}(t)$ in the trace normalization $4\pi N^{-1}\operatorname{tr}$. Thus $p_NY_{\ell0}=i q_{\ell,N}(J_3/\sqrt J)$; this phase is part of the exact Poisson--commutator identification. The finite set of nodes is $t_r=r/\sqrt J$, $r=-j,-j+1,\ldots,j$.

\begin{lemma}[Fixed-degree approximation of zonal harmonics]\label{lem:poly}
For each fixed $\ell\ge0$ and $k\ge0$, choose the degree-$\ell$ diagonal harmonic with positive leading coefficient. Then, for $N>\ell$,
\begin{equation}\label{eq:polyrate}
\norm{q_{\ell,N}-Y_{\ell0}}_{C^k([-1,1])}\le C_{\ell,k}N^{-2}.
\end{equation}
In particular, if $Q_N=\sum_{\ell=1}^d b_\ell q_{\ell,N}$ and $S_N=\sum_{\ell=1}^d -b_\ell q_{\ell,N}/[\ell(\ell+1)]$, then
\begin{equation}\label{eq:qbound}
\norm{Q_N-\omega_*}_{C^k}+\norm{S_N-\psi_*}_{C^k}\le C_{*,k}N^{-2}
\end{equation}
and their coefficients and derivatives of every fixed order remain uniformly bounded.
\end{lemma}
\begin{proof}
Index the spin nodes by $r=-j+k$, $k=0,\ldots,N-1$. Then
$y_r=2r/N=-1+(2k+1)/N$ are precisely the midpoints of an
$N$-subinterval partition of $[-1,1]$. Since
$t_r=y_r/\sqrt{1-N^{-2}}$, for every fixed nonnegative integer $s$ the composite midpoint estimate and the Taylor expansion of the scaling give
\[
\frac1N\sum_{r=-j}^{j} t_r^s=\frac12\int_{-1}^1 t^s\,dt+O_s(N^{-2}).
\]
Apply this identity to all monomials of degree at most $2\ell$. The Gram matrices of $1,t,\ldots,t^\ell$ for the discrete and continuous normalized inner products differ entrywise by $O_\ell(N^{-2})$. The limiting matrix is positive definite. The Cholesky/Gram--Schmidt map, restricted to positive-definite matrices in a fixed neighborhood of that limit, is continuously differentiable. Its coefficients therefore change by $O_\ell(N^{-2})$. The continuum orthonormal degree-$\ell$ polynomial with positive leading coefficient is $Y_{\ell0}(t)=\sqrt{(2\ell+1)/(4\pi)}P_\ell(t)$. Finite-dimensional coefficient convergence proves the same rate for all fixed derivatives on $[-1,1]$. The adjoint Casimir on a diagonal polynomial $f(J_3)$ has at spin-node $r$ the exact expression
\[
 (\mathcal C f)(r)=\bigl(J-r^2+r\bigr)\bigl(f(r)-f(r-1)\bigr)
 +\bigl(J-r^2-r\bigr)\bigl(f(r)-f(r+1)\bigr).
\]
For a degree-$\ell$ polynomial the leading coefficient of $\mathcal C f$ is $\ell(\ell+1)$ times the leading coefficient of $f$: the two second-difference terms contribute $\ell(\ell-1)$, and the first-difference terms contribute $2\ell$. The Casimir therefore preserves polynomial degree and its restriction to the polynomial flag is upper triangular with the simple diagonal eigenvalues $k(k+1)$, $0\le k\le\ell<N$. The adjoint Casimir is self-adjoint in the normalized trace inner
product, so its polynomial eigenvectors corresponding to distinct
eigenvalues are orthogonal. Since the eigenvalues along the polynomial
flag are simple, its degree-$\ell$ eigenvector is the unique
orthogonal polynomial of degree $\ell$, up to normalization and sign.
The diagonal $T^{(N)}_{\ell0}$ is the unique, up to phase, weight-zero vector in the irreducible rank-$\ell$ matrix-harmonic subspace. It is a degree-$\ell$ polynomial in $J_3$ orthogonal to all lower degrees by trace orthogonality. The choice of positive leading coefficient consequently identifies it with the preceding Gram--Schmidt polynomial. Linear combinations give~\eqref{eq:qbound}.
\end{proof}

For $m\ne0$, introduce the self-adjoint symmetric-coordinate map on matrix units $E_{rs}$ with $r-s=m$:
\begin{equation}\label{eq:BN}
B_NE_{rs}=\frac{r+s}{2\sqrt J}E_{rs},\qquad B_NZ=\frac{J_3Z+ZJ_3}{2\sqrt J}.
\end{equation}
Its domain here is the nonzero azimuthal sector. In the diagonal sector $m=0$, the anticommutator can acquire a scalar component; all coordinate estimates below are therefore formulated on the invariant nonzero-weight sectors. The linearized $m=0$ dynamics vanishes.

\begin{lemma}[Matrix-harmonic coordinate recurrence]\label{lem:rec}
With compatible Condon--Shortley phases, the exact recurrence on $1\le|m|\le\ell<N$ is
\begin{align}\label{eq:rec}
B_NT_{\ell m}^{(N)}={}&a_{\ell+1,m}\sqrt{1-\frac{\ell(\ell+2)}{N^2-1}}\,T_{\ell+1,m}^{(N)}\\[-2pt]
&+a_{\ell,m}\sqrt{1-\frac{(\ell-1)(\ell+1)}{N^2-1}}\,T_{\ell-1,m}^{(N)},\nonumber
\end{align}
where $a_{\ell,m}=\sqrt{(\ell^2-m^2)/[(2\ell-1)(2\ell+1)]}$ and unavailable terms vanish. In particular the matrix recurrence has no coupling beyond the highest degree $N-1$.
\end{lemma}
\begin{proof}
Let $j=(N-1)/2$ and $J=j(j+1)$. In the spin basis $|r\rangle$,
$-j\le r\le j$, the rank-$\ell$ matrix harmonic is
\begin{equation}\label{eq:3j}
(T_{\ell m}^{(N)})_{rs}=\sqrt{\frac{N}{4\pi}}\,
(-1)^{j-r}\sqrt{2\ell+1}
\begin{pmatrix}j&\ell&j\\-r&m&s\end{pmatrix},\qquad r-s=m.
\end{equation}
The normalization follows directly from $3j$ orthogonality:
$\sum_{r,s}(2\ell+1)|\bigl(\begin{smallmatrix}j&\ell&j\\-r&m&s\end{smallmatrix}\bigr)|^2=1$;
therefore $(4\pi/N)\operatorname{tr}[(T_{\ell m}^{(N)})^*T_{\ell m}^{(N)}]=1$.

The recurrence follows from the standard $3j$ identity
\cite[Eq.~34.3.15]{DLMF}. Cyclically permute the columns
of the $3j$ symbol in~\eqref{eq:3j} to $(\ell,j,j;m,s,-r)$.
Specializing that identity to $j_1=\ell$, $j_2=j_3=j$ and
$m_1=m$, $m_2=s$, $m_3=-r$ gives, whenever $1\le\ell<N$,
\begin{align}\label{eq:3jexplicit}
 &(2\ell+1)(r+s)
  \begin{pmatrix}\ell&j&j\\m&s&-r\end{pmatrix}\\
 &\quad=\sqrt{(N^2-\ell^2)(\ell^2-m^2)}
  \begin{pmatrix}\ell-1&j&j\\m&s&-r\end{pmatrix}\nonumber\\
 &\qquad+\sqrt{[N^2-(\ell+1)^2][(\ell+1)^2-m^2]}
  \begin{pmatrix}\ell+1&j&j\\m&s&-r\end{pmatrix}.\nonumber
\end{align}
Indeed the coefficient of the middle $3j$ in the cited recursion reduces
to $(2\ell+1)\ell(\ell+1)(r+s)$; the two coefficients on its right
contain the common factor $\ell(\ell+1)$, which has been canceled
in~\eqref{eq:3jexplicit}. Cyclic column permutations preserve the $3j$
phase, and the common factor $(-1)^{j-r}$ in~\eqref{eq:3j} cancels.
Multiply~\eqref{eq:3jexplicit} by $\sqrt{2\ell+1}$ and divide by
$2\sqrt J=\sqrt{N^2-1}$, converting the neighboring normalized $3j$
symbols with factors $\sqrt{2\ell-1}$ and $\sqrt{2\ell+3}$.
The resulting coefficients are precisely
\[
 \sqrt{\frac{\ell^2-m^2}{(2\ell-1)(2\ell+1)}}
 \sqrt{\frac{N^2-\ell^2}{N^2-1}},\qquad
 \sqrt{\frac{(\ell+1)^2-m^2}{(2\ell+1)(2\ell+3)}}
 \sqrt{\frac{N^2-(\ell+1)^2}{N^2-1}},
\]
which equal the two coefficients in~\eqref{eq:rec}.
At $\ell=|m|$ the lower coefficient is zero, and at $\ell=N-1$
the upper coefficient is zero. The identity holds for positive and
negative $m$ with the phases fixed by~\eqref{eq:3j}.
\end{proof}

Let $B_m$ denote multiplication by $x$ on $\HH_m$ and let $P_{N,m}=P_N|_{\HH_m}$. For $m\neq0$ the recurrence for $B_m$ has the same coefficients $a_{\ell,m}$ without the square-root factor.

The matrix commutator formula can now be written without implicit factors. If $Z$ is a matrix unit with $r-s=m$ and $x=(r+s)/(2\sqrt J)$, then $r/\sqrt J=x+m/(2\sqrt J)$ and $s/\sqrt J=x-m/(2\sqrt J)$. With the Poisson orientation in~\eqref{eq:lm}, the transported exact finite linearization is
\begin{align}\label{eq:exactfinite}
 A_{N,m}&=-i\mathcal G^S_{N,m}(B_{N,m})P_{N,m},
 &K_{N,m}&=-i\mathcal G^Q_{N,m}(B_{N,m})D_m^{-1}P_{N,m},\\
 \mathcal G^U_{N,m}(x)&=\sqrt J\left[U_N\!\left(x+\frac{m}{2\sqrt J}\right)-U_N\!\left(x-\frac{m}{2\sqrt J}\right)\right],
 &U_N&\in\{S_N,Q_N\}.\nonumber
\end{align}
The polynomial functional calculus in~\eqref{eq:exactfinite} is taken on the finite-dimensional $m$ block, where $B_{N,m}$ is self-adjoint. Thus $A_{N,m}$ is skew-adjoint. The inverse Laplacian is placed to the \emph{right} of the polynomial commutator multiplier; interchanging these factors would change the linearization.
\begin{lemma}[Coordinate and polynomial consistency estimates]\label{lem:coord}
On $\HH_{\ne0}:=\bigoplus_{m\ne0}\HH_m$ one has
\begin{equation}\label{eq:coord}
\norm{B_NP_N-P_NB}_{H^2\cap\HH_{\ne0}\to\HH_{\ne0}}
\le C N^{-2},\qquad B=\bigoplus_{m\ne0}B_m.
\end{equation}
For each fixed polynomial $p$,
\begin{equation}\label{eq:polyB}
\norm{p(B_N)P_N-P_Np(B)}_{H^2\cap\HH_{\ne0}\to\HH_{\ne0}}
\le C_pN^{-2}.
\end{equation}
All operators on the left are understood under the harmonic identification; $p(B_N)$ acts only within the nonzero finite azimuthal blocks.
\end{lemma}
\begin{proof}
Write $w_\ell=1+\ell(\ell+1)$ and let $E_N=B_NP_N-P_NB$, interpreted on $\HH_{\ne0}$. Fix $m\ne0$ and use the orthonormal basis $Y_{\ell m}$. All matrices are tridiagonal (and $P_N$ is diagonal). For $|m|\le\ell\le N-2$, the nonzero coefficients of $E_N$ at output degrees $\ell\pm1$ are differences of the discrete and continuous recurrence coefficients. Since
\[
0\le1-\sqrt{1-t}\le t\qquad(0\le t\le1)
\]
and $|a_{\ell,m}|\le1$, their moduli do not exceed $C\ell(\ell+2)/(N^2-1)$ or $C\ell^2/(N^2-1)$, respectively. After division by $w_\ell$ the bound is $C/N^2$, uniformly in $\ell,m,N$. The case $\ell=N-1$ has the same bound for the downward coefficient. Both upward coefficients vanish after applying $P_N$ at the output: the continuous $\ell=N$ term is discarded, and the finite coefficient vanishes already in~\eqref{eq:rec}. For an input $\ell\ge N$ only $P_NB$ may contribute; tridiagonality shows that this happens solely at $\ell=N$, and the corresponding coefficient is at most $1\le Cw_N/N^2$. 

Every column and every row of $E_N\operatorname{diag}(w_\ell^{-1})$ has at most two entries. Since each nonzero entry has modulus at most $C/N^2$, the row- and column-sum bounds and the Schur test give
\[
\sup_{m\ne0}\|E_N\operatorname{diag}(w_\ell^{-1})\|_{\HH_m\to\HH_m}\le C/N^2.
\]
This is exactly~\eqref{eq:coord}, because the spectral $H^2$ norm is the $\ell^2$ norm of $w_\ell u_{\ell m}$.

For powers, the telescoping identity is exact:
\begin{equation}\label{eq:telescope}
 B_N^kP_N-P_NB^k
 =\sum_{q=0}^{k-1}B_N^{k-1-q}(B_NP_N-P_NB)B^q,
\end{equation}
where products involving $B_N$ are restricted to the finite nonzero-$m$ sector and the equality follows by cancellation of adjacent terms. The operators $B_N$ and $B$ have $L^2$ norm at most one, while $B^q$ maps $H^2$ into itself with a norm depending only on $q$.
Indeed, multiplication by the fixed smooth function $x^q$ is bounded
on $H^2(\mathbb S^2)$ independently of the azimuthal index.
In spherical coordinates $[-\Delta_m,x]$ contains no $m^2/(1-x^2)$
term, since this potential commutes with $x$; equivalently,
the global smooth-multiplier estimate on $H^2(\mathbb S^2)$ restricts
with the same constant to every invariant azimuthal block. Apply~\eqref{eq:coord} to every summand and add the finitely many coefficients of the polynomial to prove~\eqref{eq:polyB}. The same argument controls the upper-truncation terms in the polynomial calculus.
\end{proof}

\section{Uniform azimuthal estimates and convergence of compact perturbations}\label{sec:bounds}
\begin{lemma}[Growth-rate bounds for high azimuthal modes]\label{lem:highm}
There is $M_*>0$, independent of $N>d$, such that for every $m\neq0$,
\begin{equation}\label{eq:Km}
\norm{K_m}\le\frac{\norm{\omega_*'}_\infty}{|m|+1},\qquad
\norm{K_{N,m}}\le\frac{M_*}{|m|+1}.
\end{equation}
Every eigenvalue $\lambda$ of either $L_m$ or $L_{N,m}$ satisfies $|\Rea\lambda|\le M_*/(|m|+1)$ after enlarging $M_*$ if needed. Thus for each $\delta>0$ only finitely many $m$, independently of $N$, can contribute eigenvalues with $|\Rea\lambda|\ge\delta$.
\end{lemma}
\begin{proof}
The continuous inequality follows from~\eqref{eq:lm}--\eqref{eq:D}. For the finite model, on an off-diagonal matrix unit $E_{rs}$, $r-s=m$, the quantized commutator with $Q_N(J_3/\sqrt J)$ is diagonal, with coefficient proportional to
\[
\sqrt J\,[Q_N(r/\sqrt J)-Q_N(s/\sqrt J)].
\]
The mean-value theorem bounds its absolute value by $|m|\norm{Q_N'}_\infty$. Lemma~\ref{lem:poly} makes the last factor uniform in $N$. The inverse matrix Laplacian has operator norm $[|m|(|m|+1)]^{-1}$ in the $m$ block. This gives the second bound. Finally, if $L_mu=\lambda u$ and $\norm u=1$, then
\[
|\Rea\lambda|=|\Rea\ip{K_mu}{u}|\le\norm{K_m},
\]
because the transport part is skew-adjoint. The finite model is identical. The last assertion follows by taking $|m|+1>M_*/\delta$.
\end{proof}

\begin{theorem}[Global and blockwise approximation of compact perturbations]\label{thm:K}
For every equilibrium~\eqref{eq:base} there is a constant $C_*$, independent of $N>d$, such that
\begin{equation}\label{eq:Kglobal}
\norm{K_N-K}_{\mathcal L(\HH)}\le\frac{C_*}{N}.
\end{equation}
For each fixed $m\neq0$ there is $C_{*,m}$ such that
\begin{equation}\label{eq:Kfixed}
\norm{K_{N,m}-K_m}_{\mathcal L(\HH_m)}\le\frac{C_{*,m}}{N^2}.
\end{equation}
The finite-dimensional operators are extended by zero off the harmonics of degree less than $N$.
\end{theorem}
\begin{proof}
Fix $1\le|m|<N$, and write $x=(r+s)/(2\sqrt J)$ for a matrix unit with $r-s=m$. If $\delta=m/(2\sqrt J)$, the commutator coefficient for the quantized vorticity is
\begin{equation}\label{eq:finiteDiff}
G_{N,m}(x)=\sqrt J\bigl[Q_N(x+\delta)-Q_N(x-\delta)\bigr].
\end{equation}
Taylor's formula with the odd cancellation of the centered difference yields
\[
\left|G_{N,m}(x)-mQ_N'(x)\right|
\le\frac{|m|^3}{24J}\norm{Q_N'''}_\infty
\le \frac{C_*|m|^3}{N^2}.
\]
(For degrees below three the third derivative is zero.) Inserting Lemma~\ref{lem:poly} gives, uniformly in the spectrum of $B_N$,
\begin{equation}\label{eq:Gbound}
\norm{G_{N,m}-m\omega_*'(B_N)}
\le\frac{C_*(|m|+|m|^3)}{N^2}.
\end{equation}
Here $G_{N,m}$ denotes the polynomial functional calculus of the finite self-adjoint coordinate operator $B_N$ in the given azimuthal block.

Because $D_m^{-1}$ is diagonal in harmonic degree, it commutes with $P_{N,m}$ and agrees exactly with the corresponding inverse discrete Laplacian on the retained degrees. The factor-ordering of the perturbation is essential: $D_m^{-1}$ acts first and the polynomial multiplier second. It gives
\[
K_{N,m}P_{N,m}=-iG_{N,m}D_m^{-1}P_{N,m},\qquad
P_{N,m}K_mP_{N,m}=-imP_{N,m}\omega_*'(B_m)D_m^{-1}P_{N,m}.
\]
By~\eqref{eq:Gbound} and~\eqref{eq:D},
\begin{equation}\label{eq:GxD}
\norm{[G_{N,m}-m\omega_*'(B_N)]D_m^{-1}P_{N,m}}
\le\frac{C_*(|m|+|m|^3)}{N^2|m|(|m|+1)}
\le\frac{C_*(|m|+1)}{N^2}.
\end{equation}
Moreover, the spectral Sobolev norm gives the uniform identity
\[
\norm{D_m^{-1}}_{L^2\to H^2}
=\sup_{\ell\ge|m|}\frac{1+\ell(\ell+1)}{\ell(\ell+1)}
=1+\frac1{|m|(|m|+1)}\le\frac32,\qquad m\ne0.
\]
The polynomial estimate~\eqref{eq:polyB} therefore yields
\begin{equation}\label{eq:multxD}
\norm{m\,[\omega_*'(B_N)P_{N,m}-P_{N,m}\omega_*'(B_m)]D_m^{-1}P_{N,m}}
\le\frac{C_*|m|}{N^2}.
\end{equation}
Adding~\eqref{eq:GxD} and~\eqref{eq:multxD} proves
\begin{equation}\label{eq:compressed}
\norm{K_{N,m}-P_{N,m}K_mP_{N,m}}\le \frac{C_*(|m|+1)}{N^2}.
\end{equation}

We now control the terms omitted by the harmonic truncation. Since $\omega_*'$ is a fixed polynomial of degree $d-1$, the matrix of multiplication by $\omega_*'(x)$ couples harmonics whose degrees differ by at most $d-1$. Accordingly
\begin{align*}
\norm{K_m(I-P_{N,m})}
&\le |m|\norm{\omega_*'}_\infty\sup_{\ell\ge N}\frac1{\ell(\ell+1)}
\le\frac{C_*|m|}{N^2},\\
\norm{(I-P_{N,m})K_mP_{N,m}}
&\le |m|\norm{\omega_*'}_\infty\sup_{N-d\le\ell<N}\frac1{\ell(\ell+1)}
\le\frac{C_*|m|}{N^2},
\end{align*}
where the constants can be enlarged to handle $N\le2d$. The first expression is understood as an operator on $\HH_m$ and the second exploits the finite bandwidth in degree. Thus
\begin{equation}\label{eq:fullblock}
\norm{K_{N,m}-K_m}\le\frac{C_*(|m|+1)}{N^2},\qquad 1\le|m|<N.
\end{equation}
If $|m|\ge N$, then $K_{N,m}=0$ and~\eqref{eq:Km} gives $\norm{K_{N,m}-K_m}\le C_*/N$. Taking the supremum over $m$ in the orthogonal direct sum~\eqref{eq:blocks} establishes~\eqref{eq:Kglobal}. If $m$ is fixed,~\eqref{eq:fullblock} proves~\eqref{eq:Kfixed}.
\end{proof}

\begin{corollary}[Compactness of the Euler perturbation]\label{cor:compact}
The operator $K$ is compact on $\HH$.
\end{corollary}
\begin{proof}
For every fixed $m\ne0$, $K_m=-im\omega_*'(B_m)D_m^{-1}$ is compact: $D_m^{-1}$ is diagonal with eigenvalues tending to zero, while multiplication by $\omega_*'$ is bounded. For every $M$, $\bigoplus_{|m|\le M}K_m$ is consequently compact. The first inequality of~\eqref{eq:Km} shows that these finite direct sums converge to $K$ in operator norm as $M\to\infty$. Thus $K$ is compact as a norm limit of compact operators.
\end{proof}

\section{Off-axis spectral exactness}\label{sec:spectral}
We first record the operator fact used in the proof, with both resolvent topologies stated explicitly.
\begin{lemma}[Products of compact operators and resolvents]\label{lem:compactprod}
Suppose $R_N,R$ and $R_N^*,R^*$ converge strongly, and $\sup_N(\norm{R_N}+\norm{R_N^*})<\infty$. If $C$ is compact, then
\[
\norm{(R_N-R)C}\to0,\qquad \norm{C(R_N-R)}\to0.
\]
The conclusion remains true locally uniformly in a complex parameter when the convergences and bounds hold locally uniformly.
\end{lemma}
\begin{proof}
For $C=u\otimes v^*$ of rank one,
\[
(R_N-R)C=((R_N-R)u)\otimes v^*,\qquad
C(R_N-R)=u\otimes((R_N^*-R^*)v)^*.
\]
The norms of these two operators tend to zero by the corresponding strong convergences. This proves the assertion for finite-rank $C$. Given $\varepsilon>0$, choose finite-rank $C_\varepsilon$ with $\|C-C_\varepsilon\|<\varepsilon$. Uniform boundedness of the resolvents gives
\[
\|(R_N-R)(C-C_\varepsilon)\|\le M\varepsilon,
\quad
\|(C-C_\varepsilon)(R_N-R)\|\le M\varepsilon,
\]
with $M$ independent of $N$. Letting $N\to\infty$ and then $\varepsilon\to0$ finishes the proof. On any compact parameter set the resolvent identity gives a common modulus of continuity in the parameter, since the resolvents have a common norm bound there. A finite covering then upgrades pointwise strong convergence on each fixed vector to local uniform convergence; the finite-rank approximation above preserves this uniformity.
\end{proof}

\begin{lemma}[Strong-resolvent convergence of transport operators]\label{lem:transport}
For every $z\notin i\R$,
\[
(z-A_N)^{-1}u\longrightarrow(z-A)^{-1}u\qquad (u\in\HH),
\]
and the corresponding adjoint resolvents converge strongly. The convergences are locally uniform for $z\in\C\setminus i\R$.
\end{lemma}
\begin{proof}
Let $\DD_0$ be the finite linear span of the spherical harmonics. It is a core for the skew-adjoint direct sum $A$: in each fixed $m$ block the operator $A_m$ is bounded, and truncating first in $m$ and then in harmonic degree approximates any vector in the graph norm of~\eqref{eq:Adomain}. Fix $u\in\DD_0$. It has only finitely many nonzero harmonic and
azimuthal coefficients. The recurrence~\eqref{eq:rec}, the polynomial
approximation~\eqref{eq:qbound}, and the centered difference in
\eqref{eq:exactfinite} therefore give, for all sufficiently large $N$,
\[
  \norm{(A_N-A)u}\le C_uN^{-2}.
\]
In particular $A_Nu\to Au$ in $\HH$. Both $A_N$ and $A$ are
skew-adjoint. On $(I-P_N)\HH$ the extended operator $A_N$ is zero
and its resolvent acts as multiplication by $z^{-1}$; the strong limit
uses $P_Nu\to u$ for every fixed $u\in\HH$. Consequently,
\[
\norm{(z-A_N)^{-1}},\ \norm{(z-A)^{-1}}\le |\Rea z|^{-1}.
\]
For $f=(z-A)u$, $u\in\DD_0$, the resolvent identity on the core gives
\[
(z-A_N)^{-1}f-u=(z-A_N)^{-1}(A_N-A)u\longrightarrow0.
\]
Since $(z-A)\DD_0$ is dense and the resolvents are uniformly bounded, convergence holds for every $f\in\HH$. The adjoints are $(z-A_N)^{-1*}=(\bar z+A_N)^{-1}$, and the same argument applies to $-A_N$ and $-A$. Local uniformity follows from the resolvent identity and a finite covering of compact $z$-sets.
\end{proof}

\begin{theorem}[Global off-axis spectral exactness]\label{thm:spectral}
Let $\omega_*$ be the polynomial zonal equilibrium~\eqref{eq:base}. Then:
\begin{enumerate}[label=\textnormal{(\roman*)}]
\item Every eigenvalue of $L$ in $\C\setminus i\R$ is isolated and of finite algebraic multiplicity. Every sufficiently small disk about such a value contains, for all sufficiently large $N$, the same total algebraic multiplicity of eigenvalues of $L_N$.
\item If $\lambda_{N_k}\in\spec(L_{N_k})$ belongs to a fixed compact subset of $\C\setminus i\R$ and $\lambda_{N_k}\to\lambda$, then $\lambda\in\spec(L)$.
\item On a contour $\Gamma$ whose curve and interior lie in one open half-plane and which encloses a finite off-axis spectral cluster of $L$ but no other spectrum, the Riesz projectors $\mathsf P_{N,\Gamma}$ converge in operator norm to $\mathsf P_\Gamma$.
\end{enumerate}
Consequently, no sequence of Zeitlin eigenvalues can converge to a spurious point strictly separated from the imaginary axis.
\end{theorem}
\begin{proof}
Write $R_{A,N}(z)=(z-A_N)^{-1}$ and $R_A(z)=(z-A)^{-1}$. Since $A_N,A$ are skew-adjoint, these resolvents exist for $\Rea z\ne0$. Define the analytic compact operator pencils
\[
\mathcal B_N(z)=K_NR_{A,N}(z),\qquad \mathcal B(z)=KR_A(z).
\]
By Theorem~\ref{thm:K}, $\norm{K_N-K}\to0$. Since $K$ is compact (Corollary~\ref{cor:compact}) and both $R_{A,N}$ and their adjoints converge strongly (Lemma~\ref{lem:transport}), Lemma~\ref{lem:compactprod} gives
\begin{align}\label{eq:Bconv}
\norm{\mathcal B_N(z)-\mathcal B(z)}
&\le\norm{K_N-K}\norm{R_{A,N}(z)}
  +\norm{K[R_{A,N}(z)-R_A(z)]}\longrightarrow0
\end{align}
locally uniformly off the imaginary axis. Every $\mathcal B_N(z)$ and $\mathcal B(z)$ is compact. Furthermore,
\begin{equation}\label{eq:factor}
z-L=(\id-\mathcal B(z))(z-A),\qquad
z-L_N=(\id-\mathcal B_N(z))(z-A_N).
\end{equation}
On each of the two half-planes, $\id-\mathcal B(z)$ is invertible for sufficiently large $|\Rea z|$, since $\norm{\mathcal B(z)}\le\norm K/|\Rea z|$. Analytic Fredholm theory therefore implies that its noninvertibility points are discrete, of finite algebraic type, within either half-plane. Because of~\eqref{eq:factor}, these points are precisely the off-axis eigenvalues of $L$. The same statements apply to $L_N$.

Let $\Gamma$ be a small positively oriented contour enclosing one isolated eigenvalue (or a chosen finite cluster) and lying entirely within the appropriate half-plane. On $\Gamma$, $\id-\mathcal B(z)$ has a bounded inverse. By locally uniform norm convergence in~\eqref{eq:Bconv}, the inverse $[\id-\mathcal B_N(z)]^{-1}$ exists on $\Gamma$ for all sufficiently large $N$ and converges there in operator norm to $[\id-\mathcal B(z)]^{-1}$, by the identity
\[
(\id-\mathcal B_N)^{-1}-(\id-\mathcal B)^{-1}
=(\id-\mathcal B_N)^{-1}(\mathcal B_N-\mathcal B)(\id-\mathcal B)^{-1}.
\]
In particular, the factorization~\eqref{eq:factor} and the
uniform bounds for the free resolvents yield
\[
\sup_{N\ge N_0}\sup_{z\in\Gamma}
\|(z-L_N)^{-1}\|<\infty.
\]
The same conclusion holds after restriction to any fixed azimuthal
block $\HH_m$. For norm convergence of the Riesz projectors, set
$\mathcal C_N(z)=R_{A,N}(z)K_N$, $\mathcal C(z)=R_A(z)K$. Both $\mathcal C_N\to\mathcal C$ and $\mathcal B_N\to\mathcal B$ hold in norm locally uniformly, by the two sides of Lemma~\ref{lem:compactprod}. The factorization $z-L_N=(z-A_N)(\id-\mathcal C_N(z))$ gives the resolvent correction formula
\begin{equation}\label{eq:compcorr}
(z-L_N)^{-1}-R_{A,N}(z)
= \mathcal C_N(z)[\id-\mathcal C_N(z)]^{-1}R_{A,N}(z).
\end{equation}
Indeed $(\id-\mathcal C_N)^{-1}-\id=\mathcal C_N(\id-\mathcal C_N)^{-1}$. The explicit inverse identity
\[
 (\id-\mathcal C_N)^{-1}
 =\id+R_{A,N}(\id-\mathcal B_N)^{-1}K_N
\]
shows that $\id-\mathcal C_N$ is invertible on $\Gamma$ whenever
$\id-\mathcal B_N$ is invertible.
Write $D_N(z)=\mathcal C_N(z)[\id-\mathcal C_N(z)]^{-1}$ and
$D(z)=\mathcal C(z)[\id-\mathcal C(z)]^{-1}$.
Each $D_N(z)$ and $D(z)$ is compact, and $D_N(z)\to D(z)$ in norm
locally uniformly on $\Gamma$. The decisive two-sided compactness
argument is
\[
\begin{aligned}
&\|D_N(z)R_{A,N}(z)-D(z)R_A(z)\|\\
&\qquad\le\|D_N(z)-D(z)\|\,\|R_{A,N}(z)\|
 +\|D(z)[R_{A,N}(z)-R_A(z)]\|\longrightarrow0.
\end{aligned}
\]
The second term converges in norm by Lemma~\ref{lem:compactprod}:
$D(z)$ is compact and the \emph{adjoint} transport resolvents converge
strongly. The convergence is uniform on the contour by the same lemma.
Thus the compact resolvent correction in~\eqref{eq:compcorr}
converges in operator norm, even though the free transport resolvents
need not converge in operator norm. Since the contour and its interior lie within one half-plane,
\[
\frac1{2\pi i}\int_\Gamma R_{A,N}(z)\,dz
=\frac1{2\pi i}\int_\Gamma R_A(z)\,dz=0.
\]
Integration of~\eqref{eq:compcorr} shows $\norm{\mathsf P_{N,\Gamma}-\mathsf P_\Gamma}\to0$. The azimuthal block projections commute with both $L_N$ and $L$; hence this convergence also holds on each fixed block, with algebraic multiplicities preserved blockwise even when the same eigenvalue occurs in several blocks. For large $N$, two projections at distance less than one have equal finite rank; their ranks are the total algebraic multiplicities of the eigenvalues inside the contour. This proves (i) and (iii). Finally, if $\lambda\notin\spec(L)$, norm convergence~\eqref{eq:Bconv} implies local invertibility of $\id-\mathcal B_N$ near $\lambda$ for all large $N$, ruling out a convergent spectral sequence there. This proves (ii).
\end{proof}

\begin{corollary}[Spectral counting outside a fixed growth-rate strip]\label{cor:margin}
Fix $\delta>0$ for which $L$ has no eigenvalue on the two lines $|\Rea z|=\delta$. Then the spectra of $L$ and, for all sufficiently large $N$, of $L_N$ in $\{|\Rea z|>\delta\}$ consist of finitely many eigenvalues, counted with algebraic multiplicity. The total multiplicities in the right and left sectors agree separately for all sufficiently large $N$.
\end{corollary}
\begin{proof}
Lemma~\ref{lem:highm} restricts the relevant blocks to the finite set $|m|\le M_*/\delta$. On those blocks, the transport multipliers are bounded in norm uniformly in $N$: by the mean-value theorem their commutator coefficients have modulus at most $|m|\|S_N'\|_\infty$; the fixed-degree polynomial approximation bounds these derivatives uniformly. The perturbations are bounded uniformly as well. Thus all spectra under consideration lie in a fixed disk $|z|\le R_\delta$, independently of $N$. The off-axis spectrum of $L$ in the resulting compact set is finite by Theorem~\ref{thm:spectral}. Because the two boundary lines contain no eigenvalue, they have an open neighborhood disjoint from $\sigma(L)$ inside the disk. Cover the finitely many eigenvalues with disjoint small contours and the remaining compact region with resolvent neighborhoods. Part (ii) of Theorem~\ref{thm:spectral}, applied by contradiction to a convergent subsequence, excludes finite-model eigenvalues in the remaining region. Part (i) then gives equality of the separate sums of algebraic multiplicities.
\end{proof}

\section{Quantitative approximation of isolated eigenvalues and spectral clusters}\label{sec:eig}
\begin{lemma}[Smoothness of off-axis eigenvectors]\label{lem:smooth}
If $\lambda\notin i\R$ and $L_mv=\lambda v$ for $v\in\HH_m$, then $v\in H^s(\mathbb S^2)$ for every $s\ge0$.
\end{lemma}
\begin{proof}
Rewrite the eigenvalue equation as
\[
v=(\lambda-A_m)^{-1}K_mv.
\]
The operator $D_m^{-1}$ maps $H^s\cap\HH_m$ into $H^{s+2}\cap\HH_m$. Multiplication by the fixed smooth function $\omega_*'$ preserves every Sobolev space. Since $A_m$ is multiplication by $-im\psi_*'(x)$ and $\Rea\lambda\ne0$, the inverse $(\lambda-A_m)^{-1}$ is multiplication by a smooth bounded function of $x$ (with derivatives depending on $m$ and $\lambda$), hence preserves every $H^s$. Starting with $v\in H^0$, iterate the identity to obtain $v\in H^{2k}$ for every $k$, and use interpolation for intermediate orders.
\end{proof}

\begin{lemma}[Quantitative consistency on a smooth block]\label{lem:cons}
For fixed $m\ne0$ and $v\in H^2\cap\HH_m$,
\begin{equation}\label{eq:cons}
\norm{(L_{N,m}P_{N,m}-P_{N,m}L_m)v}
\le\frac{C_{*,m}}{N^2}\norm v_{H^2}.
\end{equation}
\end{lemma}
\begin{proof}
The compact terms are bounded by the fixed-$m$ estimate in Theorem~\ref{thm:K}, together with $\norm{P_{N,m}}\le1$. For the transport terms, the commutator with the quantized stream polynomial $S_N$ is the same centered-difference operator as~\eqref{eq:finiteDiff}, with $Q_N$ replaced by $S_N$. Since $m$ is fixed, Taylor's formula bounds its difference from $m\psi_*'(B_N)$ by $C_{*,m}N^{-2}$ in $L^2$ operator norm. Lemma~\ref{lem:coord} and~\eqref{eq:qbound} then control the difference between $m\psi_*'(B_N)P_{N,m}$ and $mP_{N,m}\psi_*'(B_m)$ from $H^2$ to $L^2$ with the same rate. Adding the two parts proves~\eqref{eq:cons}. All boundary terms are included in Lemma~\ref{lem:coord}.
\end{proof}

\begin{theorem}[Second-order approximation of simple off-axis eigenvalues]\label{thm:eig}
Let $\lambda\in\spec(L)\setminus i\R$ be algebraically simple. There exists, for every sufficiently large $N$, exactly one eigenvalue $\lambda_N$ of $L_N$ in a sufficiently small fixed disk about $\lambda$, and
\begin{equation}\label{eq:eig}
|\lambda_N-\lambda|\le C_\lambda N^{-2}.
\end{equation}
\end{theorem}
\begin{proof}
An algebraically simple eigenvector $v$ lies in one nonzero azimuthal block: the projections onto the blocks commute with $L$, so the nonzero block components of $v$ are themselves eigenvectors, and algebraic simplicity excludes more than one such component. Denote its block by $m$ and normalize $\norm v=1$. Lemma~\ref{lem:smooth} gives $v\in H^2$.

By Theorem~\ref{thm:spectral}, a sufficiently small contour about $\lambda$ encloses precisely one algebraically simple eigenvalue $\lambda_N$ for large $N$. Its rank-one Riesz projector $\mathsf P_N$ converges in operator norm to $\mathsf P$, the projector associated with $\lambda$. Thus
\[
\sup_N\norm{\mathsf P_N}<\infty,\qquad
\mathsf P_NP_{N,m}v\longrightarrow\mathsf Pv=v,
\]
and consequently $\norm{\mathsf P_NP_{N,m}v}\ge1/2$ for large $N$. Since $\mathsf P_N$ commutes with $L_N$ and $L_N\mathsf P_N=\lambda_N\mathsf P_N$, we obtain
\begin{align*}
(\lambda_N-\lambda)\mathsf P_NP_{N,m}v
&=\mathsf P_N(L_N-\lambda)P_{N,m}v\\
&=\mathsf P_N(L_{N,m}P_{N,m}-P_{N,m}L_m)v.
\end{align*}
By Lemma~\ref{lem:cons}, the right side has norm at most $CN^{-2}$. Division by the established lower bound proves~\eqref{eq:eig}. The rank-one Riesz projection transfers the consistency estimate to the eigenvalue.
\end{proof}

\begin{theorem}[Eigenvalue clusters and generalized eigenspaces]\label{thm:clusters}
Let $m\ne0$ be fixed and let $\lambda\in\sigma(L_m)\setminus i\mathbb R$
be isolated, with algebraic multiplicity $\mu$ and ascent $\alpha$
(the largest Jordan-block size on its generalized eigenspace).
For all sufficiently large $N$, a fixed small disk about $\lambda$
contains exactly $\mu$ eigenvalues $\lambda_{N,1},\ldots,\lambda_{N,\mu}$
of $L_{N,m}$, counted with algebraic multiplicity, and
\begin{equation}\label{eq:cluster-rates}
\left|\frac1\mu\sum_{j=1}^{\mu}\lambda_{N,j}-\lambda\right|
 \le C_\lambda N^{-2},\qquad
\max_j|\lambda_{N,j}-\lambda|\le C_\lambda N^{-2/\alpha}.
\end{equation}
If $\mathsf P_\lambda$ and $\mathsf P_{N,\lambda}$ are the associated
Riesz projectors and $E=\operatorname{ran}\mathsf P_\lambda$, then
\begin{equation}\label{eq:subspace-rate}
\| (\mathsf P_{N,\lambda}-\mathsf P_\lambda)|_E\|
 \le C_\lambda N^{-2}.
\end{equation}
In particular, each individual eigenvalue in a semisimple off-axis
cluster converges with order $N^{-2}$.
\end{theorem}
\begin{proof}
Each generalized eigenvector in $E$ is smooth. Indeed, if
$(L_m-\lambda)v=w$ and $w\in H^s$, then
\[
 v=(\lambda-A_m)^{-1}(K_mv-w),
\]
and the regularity argument in Lemma~\ref{lem:smooth}, applied
inductively along a Jordan chain, gives $v\in H^s$ for every $s$.
As $E$ is finite-dimensional, its unit sphere has uniformly bounded
Sobolev norms of every fixed order.

Write $P=P_{N,m}$ for the harmonic truncation in this proof and
$\mathsf P_N=\mathsf P_{N,\lambda}$ for the Riesz projector. For
$v\in E$ and $z$ on the fixed spectral contour, the exact identity
\[
 (z-L_{N,m})^{-1}Pv-P(z-L_m)^{-1}v
 =(z-L_{N,m})^{-1}(L_{N,m}P-PL_m)(z-L_m)^{-1}v
\]
holds. The vector $(z-L_m)^{-1}v$ lies in $E$; its $H^2$ norm is
uniformly bounded for $z$ on the contour and $\|v\|\le1$.
The norm of the right-hand side is therefore $O(N^{-2})$ by
Lemma~\ref{lem:cons}, uniformly on that contour. Indeed,
\[
\sup_{N\ge N_0}\sup_{z\in\Gamma}
\|(z-L_{N,m})^{-1}\|<\infty
\]
follows from the Birman--Schwinger factorization and the uniform inverse
bound established in the proof of Theorem~\ref{thm:spectral}. The
blockwise preservation of algebraic multiplicity used below follows from the same theorem
restricted to $\HH_m$.
After contour integration,
\[
\|\mathsf P_NPv-Pv\|\le C_\lambda N^{-2}\|v\|.
\]
The spectral truncation satisfies
$\|(I-P)v\|\le C_\lambda N^{-2}\|v\|$ on $E$.
Because $L_{N,m}$ is extended by zero on $(I-P)\HH_m$ and the
contour is disjoint from $0$, its Riesz projector satisfies
$\mathsf P_N(I-P)=0$. Consequently $T_N=\mathsf P_NP|_E:E\to E_N:=\operatorname{ran}\mathsf P_N$
satisfies $\|T_N-I_E\|\le C_\lambda N^{-2}$, interpreted in the
ambient Hilbert space. By preservation of algebraic multiplicity, $\dim E_N=\mu$, so
$T_N$ is a uniformly bounded isomorphism for large $N$, with
uniformly bounded inverse. This proves~\eqref{eq:subspace-rate}.

Since $\mathsf P_N$ commutes with $L_{N,m}$,
\[
 L_{N,m}T_Nv-T_NL_mv
 =\mathsf P_N(L_{N,m}P-PL_m)v,
 \qquad v\in E.
\]
Thus the similarity transform of the finite restricted operator obeys
\[
 T_N^{-1}(L_{N,m}|_{E_N})T_N
 = L_m|_E+E_N^{\mathrm{err}},\qquad
 \|E_N^{\mathrm{err}}\|\le C_\lambda N^{-2}.
\]
On $E$ the operator $L_m|_E$ equals $\lambda I+J_\lambda$, where
$J_\lambda^\alpha=0$. Taking traces gives the cluster-mean estimate
in~\eqref{eq:cluster-rates}. Moreover,
$\|(z-\lambda-J_\lambda)^{-1}\|\le C|z-\lambda|^{-\alpha}$
for $0<|z-\lambda|\le1$. The eigenvalue equation under the
$O(N^{-2})$ matrix perturbation implies
$1\le C N^{-2}|\lambda_{N,j}-\lambda|^{-\alpha}$ whenever
$\lambda_{N,j}\ne\lambda$. This proves the individual estimates.
\end{proof}

\section{A degree-three unstable equilibrium}\label{sec:example}
We apply Theorem~\ref{thm:eig} to an explicit zonal equilibrium with a certified unstable eigenvalue. A scalar Schur complement yields a quantitative spectral enclosure. Let
\begin{equation}\label{eq:y30}
\omega_*(x)=Y_{30}(x)=\frac b2(5x^3-3x),\qquad b=\sqrt{\frac7{4\pi}},
\end{equation}
and fix $m=2$. The parity of the harmonic degree is invariant. In the even sector $\ell=2,4,6,\dots$, the linearization factors as
\begin{equation}\label{eq:factor30}
L_2=iGC,\qquad G=b(15B_2^2-3\id),\qquad C=\frac1{12}\id-D_2^{-1}.
\end{equation}
Here $c_2=-1/12$ whereas $c_\ell=1/12-1/[\ell(\ell+1)]>0$ for even $\ell\ge4$. Write
\[
a_\ell=\sqrt{\frac{\ell^2-4}{(2\ell-1)(2\ell+1)}},\quad
g_\ell=b[15(a_\ell^2+a_{\ell+1}^2)-3],\quad
h_\ell=15b a_{\ell+1}a_{\ell+2}.
\]
On the positive tail $\ell=4,6,\ldots$, set
$C_+=\operatorname{diag}(c_4,c_6,\ldots)$. This positive diagonal
operator is bounded and boundedly invertible, and conjugation by
$C_+^{1/2}$ gives the real symmetric Jacobi matrix $S$ with entries
\begin{equation}\label{eq:S}
S_{kk}=c_{4+2k}g_{4+2k},\qquad S_{k,k+1}=\sqrt{c_{4+2k}c_{6+2k}}\,h_{4+2k}.
\end{equation}
To obtain the scalar dispersion relation, separate the $\ell=2$ mode
from the positive tail. Directly from the coefficients above,
\[
 g_2=-\frac{6b}{7},\qquad c_2=-\frac1{12},\qquad
 c_4=\frac1{30},\qquad
 h_2=15b\,a_3a_4=\frac{5\sqrt{21}}{7\sqrt\pi}.
\]
Conjugating the tail by $C_+^{1/2}$ gives the block representation
\[
 L_2\sim
 \begin{pmatrix}
    ib/14 & i h_2\sqrt{c_4}\,e_0^* \\
    i c_2 h_2\sqrt{c_4}\,e_0 & iS
 \end{pmatrix}.
\]
Here the upper-right entry acts by
$u\mapsto i h_2\sqrt{c_4}\,\langle u,e_0\rangle$,
and the lower-left entry by
$a\mapsto i c_2h_2\sqrt{c_4}\,a e_0$.
Since $z-iS$ is invertible for $\Rea z>0$, its Schur complement is
\begin{equation}\label{eq:F}
\begin{aligned}
F(z)&=z-\frac{ib}{14}
 +c_2c_4h_2^2\ip{(z-iS)^{-1}e_0}{e_0}\\
 &=z-\frac{ib}{14}-\alpha\ip{(z-iS)^{-1}e_0}{e_0},
 \qquad \alpha=-c_2c_4h_2^2=\frac5{168\pi}.
\end{aligned}
\end{equation}
Block Gaussian elimination reduces $z-L_2$ to a direct sum of the scalar
factor $F(z)$ and the invertible factor $z-iS$, through analytic
invertible left and right multipliers in $\Rea z>0$. Hence zeros of
$F$ and eigenvalues of $L_2$ have the same algebraic multiplicities
throughout that half-plane.

\begin{proposition}[Certified unstable eigenvalue]\label{prop:cert}
For the degree-three equilibrium~\eqref{eq:y30}, the scalar function
$F$ defined by~\eqref{eq:F} has exactly one zero, counted with
multiplicity, in the disk
\begin{equation}\label{eq:disk}
|z-z_0|<10^{-6},\qquad z_0=0.07305545+0.05413579i.
\end{equation}
The corresponding eigenvalue of $L_2$ is algebraically simple and has
strictly positive real part.
\end{proposition}
\begin{proof}
Take the leading $n=100$ principal section $S_n$ and replace $(z-iS)^{-1}$ in~\eqref{eq:F} by $(z-iS_n)^{-1}$ to obtain $F_n(z)$. Since $a_\ell\le1/2$ for $m=2$ and $0<c_\ell<1/12$ on the positive tail, the off-diagonal Jacobi entries satisfy
\[
0\le S_{k,k+1}\le\beta=5b/16<0.234.
\]
Fix $\eta=0.073$ and $\gamma=0.1$. For the exponential-conjugation estimate on the infinite half-line, use the \emph{bounded} weights $W_R e_k=\exp(\gamma\min\{k,R\})e_k$. For every $R$ the difference $W_RSW_R^{-1}-S$ is tridiagonal with norm at most $2\beta(e^\gamma-1)<0.05$ by the row- and column-sum estimates. The same holds for $S_n$. Since $S$ and $S_n$ are self-adjoint, their resolvents in $\Rea z\ge\eta$ have norm at most $\eta^{-1}$. A Neumann-series perturbation therefore bounds the resolvent of each weighted operator by
\[
(\eta-0.05)^{-1}<44.
\]
Apply the argument with increasing weights and decreasing weights, then let $R\to\infty$ for any fixed pair of indices $j,k$; one obtains the uniform off-diagonal bounds
\[
|[(z-iS)^{-1}]_{jk}|,\ |[(z-iS_n)^{-1}]_{jk}|
\le44e^{-\gamma|j-k|}.
\]
Extend $S_n$ by the decoupled tail $S_{[n,\infty)}$ and denote its resolvent by $R^{\mathrm{dec}}$. Only the two interface entries $(n-1,n)$ and $(n,n-1)$ differ from $S$. For the $(0,0)$ resolvent entry the identity $R-R^{\mathrm{dec}}=R\,i(S-S^{\mathrm{dec}})R^{\mathrm{dec}}$ has only one contributing interface term, since $R^{\mathrm{dec}}_{n,0}=0$. Its modulus is at most $\beta44^2e^{-\gamma n}e^{-\gamma(n-1)}$. Multiplication by $\alpha$ proves
\begin{equation}\label{eq:tail}
|F(z)-F_n(z)|\le\alpha\beta\,44^2e^{-\gamma(2n-1)}<10^{-8}.
\end{equation}
The finite continued fraction and its derivative admit the rational interval
enclosures derived in Appendix~\ref{app:certificate}:
\begin{equation}\label{eq:intervalresults}
-2\cdot10^{-9}<\Rea F_{100}(z_0)<2\cdot10^{-9},\quad
-3\cdot10^{-9}<\Im F_{100}(z_0)<3\cdot10^{-9},\quad
\Rea F'_{100}(z_0)>1.45.
\end{equation}
In particular $|F_{100}(z_0)|<4\cdot10^{-9}$ and
$|F'_{100}(z_0)|>1.45$. The bounds use outward rational enclosures of $\pi$, square roots,
and the finite resolvent entries. The exponential factor in~\eqref{eq:tail}
is bounded using the positive rational Taylor sum
$\sum_{k=0}^{100}(19.9)^k/k!$.
Differentiating the finite resolvent twice and using $\norm{(z-iS_n)^{-1}}\le1/\eta$ yields $|F_{100}''(z)|<49$. Every point in the closed disk~\eqref{eq:disk} lies in $\Rea z>\eta$. On its boundary, Taylor's formula and~\eqref{eq:tail} show
\[
|F(z)-F'_{100}(z_0)(z-z_0)|
<4\cdot10^{-9}+\frac{49}{2}10^{-12}+10^{-8}
<1.45\cdot10^{-6}.
\]
Rouch\'e's theorem gives exactly one zero, counted with multiplicity, of $F$ in the disk; it is therefore simple. The Schur factorization proves the proposition.
\end{proof}

\paragraph{Finite-dimensional eigenvalue approximation.}
Table~\ref{tab:finite-eigs} illustrates the approximation for the
same degree-three equilibrium. The numbers are eigenvalues of the
finite $m=2$ Zeitlin blocks, evaluated from the exact polynomial
commutator formula~\eqref{eq:exactfinite}. The last column uses the
center $z_0$ of the certified continuum disk~\eqref{eq:disk}; it is a
numerical diagnostic of the second-order scaling.
\begin{table}[ht]
\centering
\caption{Zeitlin eigenvalues for $\omega_*=Y_{30}$, $m=2$.}\label{tab:finite-eigs}
\begin{tabular}{rccc}
\toprule
$N$ & $\operatorname{Re}\lambda_N$ & $\operatorname{Im}\lambda_N$
& $N^2|\lambda_N-z_0|$\\
\midrule
16 & 0.07354406 & 0.05349657 & 0.20597\\
32 & 0.07285847 & 0.05379043 & 0.40713\\
64 & 0.07299825 & 0.05405349 & 0.41053\\
128 & 0.07304116 & 0.05411524 & 0.41004\\
\bottomrule
\end{tabular}
\end{table}

\section{Consequences and further directions}\label{sec:discussion}
The results establish a quantitative correspondence between the
off-axis spectra of spherical Euler flow and its Zeitlin approximation.
For every fixed positive growth margin avoiding the Euler spectrum, the
finite models reproduce the total algebraic multiplicities of unstable and
stable eigenvalues. For isolated off-axis spectral clusters, the
Riesz-projector and eigenspace estimates recover both their algebraic
structure and quantitative approximation rates.

Allowing azimuthal indices to grow with $N$ yields the global
$O(N^{-1})$ estimate for the compact perturbation, while each fixed
nonzero block admits the sharper $O(N^{-2})$ rate. The combination of
exact coordinate recurrences, inverse-Laplacian smoothing, and compact
resolvent products provides an analytical framework for spectral
questions in other structure-preserving matrix approximations.

\appendix
\section{Rational enclosure of the scalar Schur-complement function}\label{app:certificate}
We describe a finite collection of rational inequalities that certifies
the numerical constants in Proposition~\ref{prop:cert}. Set
$n=100$, $z_0=7305545\cdot10^{-8}+i\,5413579\cdot10^{-8}$, and write
\[
 d_k=S_{kk},\qquad t_k=S_{k,k+1},\qquad 0\le k\le n-1.
\]
For $z\in\mathbb C$ define the finite continued fraction recursively by
\begin{equation}\label{eq:CF}
 r_{n-1}(z)=(z-id_{n-1})^{-1},\qquad
 r_k(z)=\bigl(z-id_k+t_k^2r_{k+1}(z)\bigr)^{-1},
 \quad k=n-2,\ldots,0.
\end{equation}
It follows from the Schur complement that
$F_n(z)=z-ib/14-\alpha r_0(z)$. The derivatives satisfy
\begin{equation}\label{eq:CFprime}
 r_{n-1}'=-r_{n-1}^2,\qquad
 r_k'=-r_k^2(1+t_k^2r_{k+1}'),\qquad
 F_n'=1-\alpha r_0'.
\end{equation}
These relations give a finite, explicit construction of the scalar
quantities used in the Rouch\'e estimate.

All coefficients are enclosed by rational intervals with directed
endpoints on the grid $10^{-70}\mathbb Z$. The enclosing interval for
$\pi$ follows from
\[
 \pi=16\arctan(1/5)-4\arctan(1/239),\qquad
 \arctan u=\sum_{j=0}^{M-1}\frac{(-1)^ju^{2j+1}}{2j+1}
       +\varepsilon_M,
 \quad |\varepsilon_M|\le\frac{|u|^{2M+1}}{2M+1},
\]
with $M=60$ and $M=22$, respectively, and the alternating-series
sign of the remainder. For any nonnegative rational interval
$[a,b]$, the square-root enclosure is obtained by integer bounds
$u^2\le a$ and $v^2\ge b$. Arithmetic operations in~\eqref{eq:CF}
and~\eqref{eq:CFprime} retain outward rational endpoints; a complex
quotient is evaluated via its conjugate and strictly positive squared
modulus. The following wider rational intervals contain the resulting
values:
\begin{align}\label{eq:rational-data}
 -1.593960259\cdot10^{-9}
 &<\Rea F_{100}(z_0)<-1.593960258\cdot10^{-9},\\
  2.582833246\cdot10^{-9}
 &<\Im F_{100}(z_0)< 2.582833247\cdot10^{-9},\\
  1.453759008705676
 &<\Rea F_{100}'(z_0)<1.453759008705677.
\end{align}
In particular,
$|F_{100}(z_0)|<4\cdot10^{-9}$ and
$|F_{100}'(z_0)|>1.45$.

For the infinite-dimensional remainder, the rational inequalities
\[
\alpha<\frac{95}{10^4},\qquad
\beta<\frac{234}{10^3},\qquad
 e^{1/10}-1<\frac{106}{10^3}
\]
show that the conjugation perturbation in the proof of
Proposition~\ref{prop:cert} has norm less than $0.05$. The positive
Taylor sum gives the exact rational lower bound
\[
 e^{199/10}>\sum_{j=0}^{100}\frac{(199/10)^j}{j!},\qquad
 \frac{(95/10^4)(234/10^3)44^2}
 {\sum_{j=0}^{100}(199/10)^j/j!}<10^{-8}.
\]
Finally,
\[
 \frac{2(95/10^4)}{(73/10^3)^3}<49,\qquad
 4\cdot10^{-9}+\frac{49}{2}\cdot10^{-12}+10^{-8}
 <1.45\cdot10^{-6}.
\]
The recurrence~\eqref{eq:CF}--\eqref{eq:CFprime}, rational bounds
\eqref{eq:rational-data}, and the displayed tail inequalities provide
the finite and infinite estimates required by Rouch\'e's theorem.

\section*{Code Availability}
The computational codes supporting the results of this study are
available from the author upon reasonable request.

\end{document}